\documentclass[a4paper,12pt]{amsart}
\usepackage{amsfonts}
\usepackage{amssymb}
\usepackage{ifthen}
\usepackage{graphicx}
\usepackage{color}
\nonstopmode \numberwithin{equation}{section}
\newtheorem{thm}{Theorem}%[section]
\newtheorem{lem}{Lemma}%[section]
\newtheorem{cor}{Corollary}[section]

\newtheorem{cl}{Claim}%[section]
\newtheorem{ca}{Case}%[section]
\newtheorem{sca}{Subcase}%[section]
\newtheorem{scl}{Subclaim}%[section]
\newtheorem{conj}{Conjecture}

\theoremstyle{definition}
\newtheorem{defn}{Definition}%[section]

\newtheorem{op}[equation]{Open Problem}
\newtheorem{ques}[equation]{Question}
\newtheorem{rem}{Remark}[section]
\newtheorem{exam}[equation]{Example}

\newcounter {own}
\def\theown {\thesection       .\arabic{own}}

\newenvironment{pf}[1][]{%
 \vskip 3mm
 \noindent
 \ifthenelse{\equal{#1}{}}%
  {{\slshape Proof. }}%
  {{\slshape #1.} }%
 }%
{\qed\bigskip}

\newcounter{alphabet}
\newcounter{tmp}
\newenvironment{Thm}[1][]{\refstepcounter{alphabet}%
\bigskip%
\noindent%
{\bf Theorem \Alph{alphabet}}%
\ifthenelse{\equal{#1}{}}{}{ (#1)}%
{\bf .} \itshape}{\vskip 8pt}

\makeatletter
\newcommand{\Ref}[1]{\@ifundefined{r@#1}{}{\setcounter{tmp}{\ref{#1}}\Alph{tmp}}}
\makeatother
\newcommand{\IB}{{\mathbb B}}

\newcommand{\IT}{{\mathbb T}}

\def\be{\begin{equation}}
\def\ee{\end{equation}}

\newcommand{\bee}{\begin{enumerate}}
\newcommand{\eee}{\end{enumerate}}

\newcommand{\blem}{\begin{lem}}
\newcommand{\elem}{\end{lem}}
\newcommand{\bthm}{\begin{thm}}
\newcommand{\ethm}{\end{thm}}
\newcommand{\bcor}{\begin{cor}}
\newcommand{\ecor}{\end{cor}}
\newcommand{\beg}{\begin{exam}}
\newcommand{\eeg}{\end{exam}}
\newcommand{\begs}{\begin{examples}}
\newcommand{\eegs}{\end{examples}}
\newcommand{\bdefe}{\begin{defn}}
\newcommand{\edefe}{\end{defn}}
\newcommand{\bprob}{\begin{prob}}
\newcommand{\eprob}{\end{prob}}
\newcommand{\bques}{\begin{ques}}
\newcommand{\eques}{\end{ques}}
\newcommand{\bei}{\begin{itemize}}
\newcommand{\eei}{\end{itemize}}
\newcommand{\bcon}{\begin{conj}}
\newcommand{\econ}{\end{conj}}
\newcommand{\bop}{\begin{op}}
\newcommand{\eop}{\end{op}}

\newcommand{\bca}{\begin{ca}}
\newcommand{\eca}{\end{ca}}
\newcommand{\bsca}{\begin{sca}}
\newcommand{\esca}{\end{sca}}

\newcommand{\bcl}{\begin{cl}}
\newcommand{\ecl}{\end{cl}}

\newcommand{\bscl}{\begin{scl}}
\newcommand{\escl}{\end{scl}}

\newcommand{\bcons}{\begin{conjs}}
\newcommand{\econs}{\end{conjs}}
\newcommand{\bprop}{\begin{propo}}
\newcommand{\eprop}{\end{propo}}
\newcommand{\br}{\begin{rem}}
\newcommand{\er}{\end{rem}}
\newcommand{\brs}{\begin{rems}}
\newcommand{\ers}{\end{rems}}
\newcommand{\bo}{\begin{obser}}
\newcommand{\eo}{\end{obser}}
\newcommand{\bos}{\begin{obsers}}
\newcommand{\eos}{\end{obsers}}
\newcommand{\bpf}{\begin{pf}}
\newcommand{\epf}{\end{pf}}
\newcommand{\ba}{\begin{array}}
\newcommand{\ea}{\end{array}}
\newcommand{\beq}{\begin{eqnarray}}
\newcommand{\beqq}{\begin{eqnarray*}}
\newcommand{\eeq}{\end{eqnarray}}
\newcommand{\eeqq}{\end{eqnarray*}}

\newcommand{\ra}{\rightarrow}

\newcommand{\ds}{\displaystyle}

\newcounter{minutes}
\divide\time by 60
\newcounter{hours}
\multiply\time by 60 \addtocounter{minutes}{-\time}
\begin{document}

\bibliographystyle{amsplain}
\title{ Schwarz-Pick type estimates of pluriharmonic mappings in the unit polydisk}

%%%%%%%% BEGIN TIMESTAMP
\def\thefootnote{}
\footnotetext{ \texttt{\tiny File:~\jobname .tex,
          printed: \number\day-\number\month-\number\year,
          \thehours.\ifnum\theminutes<10{0}\fi\theminutes}
} \makeatletter\def\thefootnote{\@arabic\c@footnote}\makeatother
%%%%%%%% END TIMESTAMP

\author{Shaolin Chen}
\address{Shaolin Chen, Department of Mathematics and Computational
Science, Hengyang Normal University, Hengyang, Hunan 421008,
People's Republic of China.} \email{mathechen@126.com}

\author{Antti Rasila}
\address{Antti Rasila,
Department of Mathematics and Systems Analysis, Aalto University, P.
O. Box 11100,  FI-00076 Aalto,  Finland.
 }
 \email{antti.rasila@iki.fi}

%\author{Matti Vuorinen}
%\address{Matti Vuorinen,
%Department of Mathematics,
% University of Turku, Turku 20014, Finland.
% }
 %\email{vuorinen@utu.fi}

%\author{S. Ponnusamy $^\dagger $
%${}^{~\mathbf{*}}$}
%\address{S. Ponnusamy,
%Department of Mathematics,
%Indian Institute of Technology Madras, Chennai-600 036, India.}
%\email{samy@iitm.ac.in}

%\author{X. Wang $^\dagger {}^\dagger$
%$^{\mathbf{*}}$
%${}^{~\mathbf{*}}$
%}
%\address{X. Wang, Department of Mathematics,
%Hunan Normal University, Changsha, Hunan 410081, People's Republic
%of China.} \email{xtwang@hunnu.edu.cn}

\keywords{Pluriharmonic   mapping, Schwarz-Pick type estimate,  polydisk.\\
%${}^{\mathbf{*}}
%^\dagger$ {\tt  This author is currently at the
%Indian Statistical Institute (ISI), Chennai Centre, SETS (Society
%for Electronic Transactions and security), MGR Knowledge City, CIT
%Campus, Taramani, Chennai 600 113, India. %\email{samy@isichennai.res.in}
%}\\
%$^\dagger {}^\dagger %{}^{\mathbf{*}}
%$ Corresponding author
}
%\date{\today } %November 4, 10;
%File: Ch-W-S12${}_{}$equiv-mod${}_{}$submit.tex}

\begin{abstract}
In this paper, we will give  Schwarz-Pick type estimates of
arbitrary order partial derivatives for bounded pluriharmonic
mappings defined in the unit polydisk. Our main results are
generalizations of results of Colonna for planar harmonic mappings
in [Indiana Univ. Math. J. 38: 829--840, 1989].

\end{abstract}
%Mestrovi\'c}, An isoperimetric type inequality for harmonic
%functions, \textit{J. Math. Anal. Appl.,} {\bf 373}(2011), 439--448.

%\thanks{The research was partly supported by
%NSF of China (No. 11071063)} %and  Hunan Provincial Innovation
%%Foundation for Postgraduate (No. 125000-4113).  }

\maketitle \pagestyle{myheadings} \markboth{ Sh. Chen and A.
Rasila}{ Schwarz-Pick type estimates of pluriharmonic mappings in
the unit polydisk}

\section{Introduction and main results  }\label{csw-sec1}
Let $\mathbb{C}^{n}$ denote the complex Euclidean $n$-space. For
$z=(z_{1},\ldots,z_{n})\in \mathbb{C}^{n}$, the conjugate of $z$,
denoted by $\overline{z}$, is defined by
$\overline{z}=(\overline{z}_{1},\ldots, \overline{z}_{n} ).$ %where
%$T$ is the matrix transpose.
 For $z$ and
$w=(w_{1},\ldots,w_{n})\in\mathbb{C}^{n} $, the  {\it inner product}
on $\mathbb{C}^{n}$ and the {\it Euclidean norm} of $z$ are given by
$\langle z,w\rangle := \sum_{k=1}^nz_k\overline{w}_k$ and $
\|z\|:={\langle z,z\rangle}^{1/2}, $ respectively. For $a\in
\mathbb{C}^n$, $\IB^n(a,r)=\{z\in \mathbb{C}^{n}:\, \|z-a\|<r\} $ is
the (open) ball of radius $r$ with center $a$. Also, we let
$\IB^n(r):=\IB^n(0,r)$ and use $\IB^n$ to denote the unit ball
$\IB^n(1)$, and $\mathbb{D}=\mathbb{B}^1$. Let
$\mathbb{D}^{n}=\mathbb{D}\times\cdots\times\mathbb{D}~\mbox{($n$
times)}$ be the polydisc in $\mathbb{C}^{n}$ and $\IT^{n}=\IT\times
\cdots\times \IT~\mbox{($n$ times)}$, where $\IT$ is the unit circle
in $\mathbb{C}^{1}$. A multi-index $k=(k_{1},\ldots,k_{n})$ consists
of $n$ nonnegative integers $k_{j}$, where $j\in\{1,\ldots,n\}$. The
degree of a multi-index $k$ is the sum $|k|=\sum_{j=1}^{n}k_{j}$.
Given another multi-index $\alpha=(\alpha_{1},\ldots,\alpha_{n})$,
let $k^{\alpha}=(k_{1}^{\alpha_{1}},\ldots,k_{n}^{\alpha_{n}})$. For
$z=(z_{1},\ldots,z_{n})\in \mathbb{C}^{n}$, let
$\|z\|=\left(\sum_{j=1}^{n}|z_{j}|^{2}\right)^{1/2}$,
$\|z\|_{\infty}=\max_{1\leq j\leq n}|z_{j}|$ and
$z^{k}=\Pi_{k=1}^{n}z_{j}^{k_{j}}.$

A continuous complex-valued function $f$ defined on a domain
$\Omega\subset\mathbb{C}^{n}$ is said to be {\it pluriharmonic} if
for each fixed $z\in \Omega$ and $\theta\in\partial\mathbb{B}^{n}$,
the function $f(z+\theta\zeta)$ is harmonic in $\{\zeta:\;
\|\zeta\|< d_{\Omega}(z)\}$, where $d_{\Omega}(z)$ denotes the
distance from $z$ to the
boundary $\partial\Omega$ of $\Omega$ (cf. \cite{R2}). %It follows from
%\cite[Theorem 4.4.9]{R} that a real-valued function $u$ defined on
%$G$ is pluriharmonic if and only if $u$ is the real part of a
%holomorphic function on $G$.
%Clearly,
If $\Omega\subset\mathbb{C}^{n}$ is a simply connected domain, then
a function $f:\,\Omega\ra \mathbb{C}$ is pluriharmonic if and only
if $f$ has a representation $f=h+\overline{g},$ where  $h$ and $g$
are holomorphic in $\Omega$ (see \cite{Vl}). Let
$\mathcal{P}(\Omega,\mathbb{C}^{N})$ be the class of all
pluriharmonic mappings $f=(f_{1},\ldots, f_{N})$ from a domain
$\Omega\subset\mathbb{C}^{n}$ to $\mathbb{C}^{N}$, where $N$ is a
positive integer and $f_{j}~(1\leq j\leq N)$ are pluriharmonic
mappings from $\Omega$ into $\mathbb{C}.$ For a mapping
$f\in\mathcal{P}(\Omega,\mathbb{C}^{N}),$ we use $Df$ and
$\overline{D}f$ to denote the two $N\times n$ matrices
$\left(\partial f_{j}/\partial z_{m}\right)_{N\times n}$ and
$\left(\partial f_{j}/\partial \overline{z}_{m}\right)_{N\times n}$,
respectively.
 We refer to \cite{HG1,CPW-1,CPW-2,CPW-3,DHK,I} for more details on
pluriharmonic mappings. In particular, if $n=1$, then pluriharmonic
mappings are planar harmonic mappings (cf.
\cite{Clunie-Small-84,Du}). Therefore, pluriharmonic mappings can be
understood as the natural generalization of planar harmonic mappings
to several complex variables.

We first recall  the classical Schwarz Lemma for analytic functions
$f$ of $\mathbb{D}$ into itself: \be\label{eq-CPR1}
|f'(z)|\leq\frac{1-|f(z)|^{2}}{1-|z|^{2}},~z\in\mathbb{D}. \ee In
1920, Sz\'asz \cite{S} extended the inequality (\ref{eq-CPR1}) to
the following estimate involving higher order derivatives:

\be\label{eq-CPR2}
|f^{(2m+1)}(z)|\leq\frac{(2m+1)!}{(1-|z|^{2})^{2m+1}}\sum_{k=0}^{m}{m\choose
k}^{2}|z|^{2k}, \ee where $m\in\{1,2,\ldots\}.$ In 1985, Ruscheweyh
(cf. \cite{ AW,AW-09, R}) improved (\ref{eq-CPR2}) to the following
sharp estimate: \be\label{eq-CPR3}
|f^{(n)}(z)|\leq\frac{n!(1-|f(z)|^{2})}{(1-|z|)^{n}(1+|z|)}. \ee
Recently, the inequality (\ref{eq-CPR3}) was  generalized into a
variety of forms (see \cite{ADR,AR,AW-09,DCP,DP,LC,WL}).

In 1989, Colonna established an analogue of the Schwarz-Pick lemma
for planar harmonic mappings, which is the following.

\begin{Thm}{\rm (\cite[Theorems 3 \mbox{and} 4]{Co})}\label{Co}
Let $f$ be a harmonic mapping of $\mathbb{D}$ into $\mathbb{D}$.
Then for $z\in\mathbb{D}$,
$$\left|\frac{\partial f(z)}{\partial z}\right|+\left|\frac{\partial f(z)}{\partial\overline{z}}\right|\leq\frac{4}{\pi}\frac{1}{1-|z|^{2}}.
$$
This estimate is sharp, and all the extremal functions are
$$f(z)=\frac{2\gamma }{\pi}\arg \left (
\frac{1+\psi(z)}{1-\psi(z)}\right), $$ where $|\gamma|=1$ and $\psi$
is a conformal automorphism of $\mathbb{D}$.
\end{Thm}
We refer to \cite{HG1,H1,H2,CPW0,CPW-BMMSC2011,CPW-1,CPRW,KV,LLZ,SW}
for further discussion on this topic.

In this paper, we   generalize Theorem \Ref{Co} to higher
dimensional
case, and give the estimate for the partial derivatives of arbitrary order. %Since the high
%dimensional case and the partial derivatives of arbitrary order is
 %different from the one dimensional case, the methods of the
%proof used in \cite{Co} are can't be used.
One should note that the higher dimensional case is very different
from the one dimensional situation and,  because we are dealing with
partial derivatives of arbitrary order, the method of proof from
\cite{Co} can not be used. By using the coefficient estimates and
the Cauchy integral formula, we prove the following result.

  %Our result is
%the following.

\begin{thm}\label{thm-1}
Let $f\in\mathcal{P}(\mathbb{D}^{n},\mathbb{D}).$ Then
$$\left|\frac{\partial^{\alpha} f(z)}{\partial
z_{1}^{\alpha_{1}}\cdots
\partial z_{n}^{\alpha_{n}}}\right|+\left| \frac{\partial^{\alpha} f(z)}{\partial
\overline{z}_{1}^{\alpha_{1}}\cdots
\partial \overline{z}_{n}^{\alpha_{n}}}\right|\leq\alpha!\frac{4}{\pi}\frac{(1+\|z\|_{\infty})^{|\alpha|-n}}{(1-\|z\|_{\infty}^{2})^{|\alpha|}},$$
where  $\alpha=(\alpha_{1},\cdots,\alpha_{n})$ is a multi-index with
$\alpha_{j}>0,~j\in\{1,\ldots,n\}$.
\end{thm}
We remark that if $|\alpha|=n=1$, then Theorem \ref{thm-1} coincides
with Theorem \Ref{Co}.

%In 1959, Heinz prove the Schwarz lemma of  planar harmonic mappings
%as follows.

%\begin{Thm}{\rm (\cite[Lemma]{HZ})}\label{H}
%Let $f$ be a harmonic mapping of $\mathbb{D}$ into $\mathbb{D}$ with
%$f(0)=0$. Then for $z\in\mathbb{D}$,
%$$|f(z)|\leq\frac{4}{\pi}\arctan|z|.
%$$
%\end{Thm}

It is well-known  that there are no biholomorphic mappings between
$\mathbb{D}^{n}$ and $\mathbb{B}^{n}$ (cf. \cite{R1,R2}). Hence
pluriharmonic mappings between $\mathbb{D}^{n}$ and $\mathbb{B}^{n}$
are particularly  interesting in the theory of several complex
variables. The following result is an analogue of \cite[Theorem
4]{HG1} for vector-valued pluriharmonic mappings defined in
$\mathbb{B}^{n}$.

\begin{thm}%\label{thm-2}  
If $f\in\mathcal{P}(\mathbb{D}^{n},\mathbb{B}^{N})$, then

\be\label{eq17}\max_{\theta\in\mathbb{C}^{n},~\|\theta\|_{\infty}=1}\left\|Df(z)\theta+
\overline{D}f(z)\overline{\theta}\right\|\leq\frac{4}{\pi(1-\|z\|_{\infty}^{2})},\ee
where $\theta$ is regarded as a column vector.\\

 If
$f\in\mathcal{P}(\mathbb{D}^{n},\mathbb{B}^{N})$ with $f(0)=0$, then
\be\label{eq18} \|f(z)\|\leq\frac{4}{\pi}\arctan\|z\|_{\infty}.\ee

%\item{{\rm(2)}}~If $f\in\mathcal{P}(\mathbb{B}^{n},\mathbb{B}^{N})),$
%then
%$$\max_{\theta\in\partial\mathbb{B}^{n}}\left|Df(z)\theta+\overline{D}f(z)\overline{\theta}\right|\leq\frac{4}{\pi(1-|z|^{2})}.$$
%where $\theta\in\partial\mathbb{B}^{n}.$
\end{thm}

We remark that if $n=N=1$, then the estimates  %Theorem \ref{thm-2}
(\ref{eq17}) and
  (\ref{eq18}) coincide with
Theorem \Ref{Co} and \cite[Lemma]{HZ}, respectively. %In fact, Theorem \ref{thm-2} (2) was prove in
%\cite{HG1} by using a different proof method.

%The following result is a improvement of \cite[Lemma 1]{QW}.

%\begin{thm}\label{thm-4}
%Let $u$ be a pluriharmonic from $\mathbb{B}^{n}$ into $\mathbb{R}$
%with $|u|<1$. Then for each
%$\theta=(\theta_{1},\cdots,\theta_{n})^{T}\in\partial\mathbb{B}^{n}$,
%we have
%$$\left|\sum_{k=1}^{n}\frac{\partial u(z)}{\partial
%z_{k}}\theta_{k}\right|+\left|\sum_{k=1}^{n}\frac{\partial
%u(z)}{\partial
%\overline{z}_{k}}\overline{\theta}_{k}\right|\leq\frac{4\left(1-u^{2}(z)\right)}{\pi(1-|z|^{2})}.$$

%\end{thm}

\section{The proofs of the main results  }\label{csw-sec2}

\begin{lem}\label{lem-1}
Let $m$ be a positive integer and $\gamma$ be a real constant. Then
\[
\int_{0}^{2\pi}|\cos (m\theta+\gamma)|d\theta=4.
\]
\end{lem}
\bpf By elementary calculations, we have

\begin{eqnarray*}
\int_{0}^{2\pi}|\cos
(m\theta+\gamma)|d\theta&=&\frac{1}{m}\int_{\gamma}^{2m\pi+\gamma}|\cos t|\,dt\\
&=&\frac{1}{m}\sum_{k=1}^{2m}\int_{(k-1)\pi+\gamma}^{k\pi+\gamma}|\cos t|\,dt\\
&=&\frac{1}{m}\sum_{k=1}^{2m}\int_{0}^{\pi}|\cos t|\,dt\\
&=&2\int_{0}^{\pi}|\cos t|\,dt\\
&=&4.
\end{eqnarray*}
The proof of the lemma is complete. \epf

\subsection*{Proof of Theorem \ref{thm-1}} Since $\mathbb{D}^{n}$ is a simply connected domain in
$\mathbb{C}^{n}$, we see that $f$  has a representation
$f=h+\overline{g}$, where  $h$ and $ g$ are holomorphic in
$\mathbb{D}^{n}$. Let $k=(k_{1},\ldots,k_{n})$ be a multi-index.
   Then $f$ can be expressed as a power series as follows
$$f(z)=h(z)+\overline{g(z)}=
\sum_{k}a_{k}z^{k}+\sum_{k}\overline{b}_{k}\overline{z}^{k}.$$
%where
%$f_{j}(z)=f_{j}(0)+\sum_{\alpha}a_{j,\alpha}z^{\alpha}+\sum_{\alpha}\overline{b}_{j,\alpha}\overline{z}^{\alpha}~(1\leq
%j\leq N)$, $a_{\alpha}=(a_{1,\alpha},\cdots, a_{N,\alpha})$ and
%$b_{\alpha}=(b_{1,\alpha},\cdots, b_{N,\alpha}).$

%\vspace{7pt}
\bcl  {\em For  $|k|\geq1$,
$|a_{k}|+|b_{k}|\leq\frac{4}{\pi}.$} \ecl

Now we prove Claim 1. Let
$z=(z_{1},\ldots,z_{n})=(r_{1}e^{i\theta_{1}},\ldots,r_{n}e^{i\theta_{n}})\in
\mathbb{D}^{n}$, where  $0\leq r_{j}<1$ for all
$j\in\{1,\ldots,n\}.$ Then for $|k|\geq1$,

 \be\label{eq-1}a_{k}r_{1}^{k_{1}}\cdots
r_{n}^{k_{n}}=\frac{1}{(2\pi)^{n}}\int_{0}^{2\pi}\cdots
\int_{0}^{2\pi}f(r_{1}^{k_{1}}e^{i\theta_{1}},\ldots,r_{n}^{k_{n}}e^{i\theta_{n}})
e^{-i\sum_{j=1}^{n}k_{j}\theta_{j}}d\theta_{1}\cdots d\theta_{n}\ee
and

 \be\label{eq-2}\overline{b}_{k}r_{1}^{k_{1}}\cdots
r_{n}^{k_{n}}=\frac{1}{(2\pi)^{n}}\int_{0}^{2\pi}\cdots
\int_{0}^{2\pi}f(r_{1}^{k_{1}}e^{i\theta_{1}},\ldots,r_{n}^{k_{n}}e^{i\theta_{n}})
e^{i\sum_{j=1}^{n}k_{j}\theta_{j}}d\theta_{1}\cdots d\theta_{n}.\ee
By (\ref{eq-1}) and (\ref{eq-2}), we get

%\begin{eqnarray*}
\beq\label{eq-3}  &&r_{1}^{k_{1}}\cdots
r_{n}^{k_{n}}(|a_{k}|+|b_{k}|)\\ \nonumber
&=&\bigg|\frac{1}{(2\pi)^{n}}\int_{0}^{2\pi}\cdots
\int_{0}^{2\pi}\big(e^{-i\sum_{j=1}^{n}k_{j}\theta_{j}}
e^{-i\arg a_{k}}\\
\nonumber &&+e^{i\sum_{j=1}^{n}k_{j}\theta_{j}} e^{i\arg
b_{k}}\big)f(r_{1}^{k_{1}}e^{i\theta_{1}},\ldots,r_{n}^{k_{n}}e^{i\theta_{n}})d\theta_{1}\cdots
d\theta_{n}\bigg|\\
\nonumber&\leq&\frac{1}{(2\pi)^{n}}\int_{0}^{2\pi}\cdots
\int_{0}^{2\pi}\big|1+e^{(2\sum_{j=1}^{n}k_{j}\theta_{j}+\arg
a_{k}+\arg b_{k})i}\big|\\ \nonumber
&&\times\big|f(r_{1}^{k_{1}}e^{i\theta_{1}},\ldots,r_{n}^{k_{n}}e^{i\theta_{n}})\big|d\theta_{1}\cdots
d\theta_{n}\\
\nonumber&\leq&\frac{2}{(2\pi)^{n}}\int_{0}^{2\pi}\cdots
\int_{0}^{2\pi}\left|\cos\Big[\sum_{j=1}^{n}k_{j}\theta_{j}+\frac{(\arg
a_{k}+\arg b_{k})}{2}\Big]\right|d\theta_{1}\cdots d\theta_{n}. \eeq
%\end{eqnarray*}
Since $|k|\geq1$, without loss of generality, we  assume that
$k_{1}\neq0$. By
 using
Lemma \ref{lem-1}, we see that

\be\label{eq-4}
\int_{0}^{2\pi}\Bigg|\cos\left[\sum_{j=1}^{n}k_{j}\theta_{j}+\frac{(\arg
a_{k}+\arg b_{k})}{2}\right]\Bigg|d\theta_{1}=4.\ee Then
(\ref{eq-3}) and (\ref{eq-4}) yield that

$$r_{1}^{k_{1}}\cdots
r_{n}^{k_{n}}(|a_{k}|+|b_{k}|)\leq\frac{4}{\pi}.$$ For
$j\in\{1,\ldots,n\}$, by letting $r_{j}\rightarrow1-$, we obtain the
desired result. %\be\label{eq-5}|a_{k}|+|b_{k}|\leq\frac{4}{\pi}.\ee%

For $j\in\{1,\ldots,n\}$ and
$z=(z_{1},\ldots,z_{n})\in\mathbb{D}^{n},$ let
$$\phi(\zeta)=(\phi_{1}(\zeta_{1}),\ldots,\phi_{n}(\zeta_{n})),$$
where $\zeta=(\zeta_{1},\ldots,\zeta_{n})\in\mathbb{D}^{n}$ and
$$\phi_{j}(\zeta_{j})=\frac{z_{j}+\zeta_{j}}{1+\overline{z_{j}}\zeta_{j}}.$$
Then $f\circ\phi$ can be written as following form
$$T(\zeta)=f(\phi(\zeta))=H(\zeta)+\overline{G(\zeta)}=\sum_{k}c_{k}\zeta^{k}+\sum_{k}\overline{d}_{k}\overline{\zeta}^{k},$$
where $H=h\circ\phi$ and  $G=g\circ\phi$. By using the  proof of
Claim 1, we get \be\label{eq-9} |c_{k}|+|d_{k}|\leq\frac{4}{\pi}.\ee
 For $r\in(0,1)$ and
$z\in\mathbb{D}^{n}$ with $\|z\|_{\infty}<r$, by the Cauchy integral
formula (cf. \cite{R1,WL}), we see that

\beq\label{eq-6} f(z)&=&\frac{1}{(2\pi
i)^{n}}\int_{|\eta_{1}|=r}\cdots\int_{|\eta_{n}|=r}
\frac{h(\eta_{1},\ldots,\eta_{n})}{\prod_{j=1}^{n}(\eta_{j}-z_{j})}d\eta_{1}\cdots
d\eta_{n}\\ \nonumber &&+\overline{\frac{1}{(2\pi
i)^{n}}\int_{|\eta_{1}|=r}\cdots\int_{|\eta_{n}|=r}
\frac{g(\eta_{1},\ldots,\eta_{n})}{\prod_{j=1}^{n}(\eta_{j}-z_{j})}d\eta_{1}\cdots
d\eta_{n}},\eeq which implies that

\be\label{eq-7} \frac{\partial^{\alpha} f(z)}{\partial
z_{1}^{\alpha_{1}}\cdots
\partial z_{n}^{\alpha_{n}}}=\frac{\alpha!}{(2\pi
i)^{n}}\int_{|\eta_{1}|=r}\cdots\int_{|\eta_{n}|=r}\frac{h(\eta_{1},\ldots,\eta_{n})}{\prod_{j=1}^{n}(\eta_{j}-z_{j})^{\alpha_{j}+1}}d\eta_{1}\cdots
d\eta_{n}\ee and

\be\label{eq-8} \frac{\partial^{\alpha} f(z)}{\partial
\overline{z}_{1}^{\alpha_{1}}\cdots
\partial\overline{ z}_{n}^{\alpha_{n}}}=\overline{\frac{\alpha!}{(2\pi
i)^{n}}\int_{|\eta_{1}|=r}\cdots\int_{|\eta_{n}|=r}\frac{g(\eta_{1},\ldots,\eta_{n})}{\prod_{j=1}^{n}(\eta_{j}-z_{j})^{\alpha_{j}+1}}d\eta_{1}\cdots
d\eta_{n}}.\ee For $j\in\{1,\ldots,n\}$, by taking
$\eta_{j}=\phi_{j}(\zeta_{j})=\frac{z_{j}+\zeta_{j}}{1+\overline{z_{j}}\zeta_{j}}$, we see from
(\ref{eq-9}), (\ref{eq-7}) and (\ref{eq-8}) that

\begin{eqnarray*}
 \frac{\partial^{\alpha} f(z)}{\partial
z_{1}^{\alpha_{1}}\cdots
\partial z_{n}^{\alpha_{n}}}&=&\frac{\alpha!}{(2\pi
i)^{n}\prod_{j=1}^{n}(1-|z_{j}|^{2})^{\alpha_{j}}}\\
&&\times\int_{|\phi_{1}(\zeta_{1})|=r}\cdots\int_{|\phi_{n}(\zeta_{n})|=r}
\frac{H(\zeta_{1},\ldots,\zeta_{n})\prod_{j=1}^{n}(1+\overline{z}_{j}\zeta_{j})^{\alpha_{j}-1}}{\prod_{j=1}^{n}\zeta_{j}^{\alpha_{j}+1}}d\zeta_{1}\cdots
d\zeta_{n}\\
&=&\frac{\alpha!}{\prod_{j=1}^{n}(1-|z_{j}|^{2})^{\alpha_{j}}}\\
&&\times\sum_{k_{1}=0}^{\alpha_{1}-1}\cdots\sum_{k_{n}=0}^{\alpha_{n}-1}{\alpha_{1}-1\choose
k_{1}}\cdots{\alpha_{n}-1\choose k_{n}}
c_{\alpha_{1}-k_{1},\ldots,\alpha_{n}-k_{n}}\overline{\prod_{j=1}^{n}z_{j}^{k_{j}}}
\end{eqnarray*}
and
\begin{eqnarray*}
 \frac{\partial^{\alpha} f(z)}{\partial
\overline{z}_{1}^{\alpha_{1}}\cdots
\partial \overline{z}_{n}^{\alpha_{n}}}&=&\frac{\alpha!}{\overline{(2\pi
i)^{n}}\prod_{j=1}^{n}(1-|z_{j}|^{2})^{\alpha_{j}}}\\
&&\times\overline{\int_{|\phi_{1}(\zeta_{1})|=r}\cdots\int_{|\phi_{n}(\zeta_{n})|=r}
\frac{G(\zeta_{1},\ldots,\zeta_{n})\prod_{j=1}^{n}(1+\overline{z}_{j}\zeta_{j})^{\alpha_{j}-1}}{\prod_{j=1}^{n}\zeta_{j}^{\alpha_{j}+1}}d\zeta_{1}\cdots
d\zeta_{n}}\\
&=&\frac{\alpha!}{\Pi_{j=1}^{n}(1-|z_{j}|^{2})^{\alpha_{j}}}\\
&&\times\sum_{k_{1}=0}^{\alpha_{1}-1}\cdots\sum_{k_{n}=0}^{\alpha_{n}-1}{\alpha_{1}-1\choose
k_{1}}\cdots{\alpha_{n}-1\choose k_{n}}
\overline{d}_{\alpha_{1}-k_{1},\ldots,\alpha_{n}-k_{n}}\prod_{j=1}^{n}z_{j}^{k_{j}}.
\end{eqnarray*} It follows that

\begin{eqnarray*}
 &&\left|\frac{\partial^{\alpha} f(z)}{\partial
z_{1}^{\alpha_{1}}\cdots
\partial z_{n}^{\alpha_{n}}}\right|+\left| \frac{\partial^{\alpha} f(z)}{\partial
\overline{z}_{1}^{\alpha_{1}}\cdots
\partial \overline{z}_{n}^{\alpha_{n}}}\right|\\
&\leq&\frac{\alpha!}{\prod_{j=1}^{n}(1-|z_{j}|^{2})^{\alpha_{j}}}
\sum_{k_{1}=0}^{\alpha_{1}-1}\cdots\sum_{k_{n}=0}^{\alpha_{n}-1}{\alpha_{1}-1\choose
k_{1}}\cdots{\alpha_{n}-1\choose k_{n}}\\
&&\times
\big(|c_{\alpha_{1}-k_{1},\cdots,\alpha_{n}-k_{n}}|+|d_{\alpha_{1}-k_{1},\cdots,\alpha_{n}-k_{n}}|\big)\prod_{j=1}^{n}|z_{j}|^{k_{j}}\\
&\leq&\frac{4}{\pi}\frac{\alpha!}{\prod_{j=1}^{n}(1-|z_{j}|^{2})^{\alpha_{j}}}
\sum_{k_{1}=0}^{\alpha_{1}-1}\cdots\sum_{k_{n}=0}^{\alpha_{n}-1}{\alpha_{1}-1\choose
k_{1}}\cdots{\alpha_{n}-1\choose
k_{n}}\prod_{j=1}^{n}|z_{j}|^{k_{j}}\\
&\leq&\frac{4}{\pi}\frac{\alpha!}{\prod_{j=1}^{n}(1-|z_{j}|^{2})^{\alpha_{j}}}\prod_{j=1}^{n}(1+|z_{j}|)^{\alpha_{j}-1}\\
&=&\alpha!\frac{4}{\pi}\prod_{j=1}^{n}\frac{(1+|z_{j}|)^{\alpha_{j}-1}}{(1-|z_{j}|^{2})^{\alpha_{j}}}\\
&\leq&\alpha!\frac{4}{\pi}\frac{(1+\|z\|_{\infty})^{|\alpha|-n}}{(1-\|z\|_{\infty}^{2})^{|\alpha|}}.
\end{eqnarray*}
The proof of the theorem is complete. \hfill $\Box$

%~{\rm(\mbox{or}}~z\in\mathbb{B}^{n})

\begin{lem}\label{thm-2}
Let $\beta=(\beta_{1},\ldots,\beta_{n})$ be a multi-index consisting
of $n$ nonnegative integers $\beta_{v}$ and
$f=(f_{1},\ldots,f_{N})\in\mathcal{P}(\mathbb{D}^{n},\mathbb{B}^{N}),$
where $v\in\{1,\ldots,n\}$. Suppose that for $z\in\mathbb{D}^{n},$
$f(z)=\sum_{\beta}a_{\beta}z^{\beta}+\sum_{\beta}\overline{b}_{\beta}\overline{z}^{\beta}$,
where  for $j\in\{1,\ldots,N\}$,
$f_{j}=\sum_{\beta}a_{j,\beta}z^{\beta}+\sum_{\beta}\overline{b}_{j,\beta}\overline{z}^{\beta}$,
$a_{\beta}=(a_{1,\beta},\ldots,a_{N,\beta})$ and
$b_{\beta}=(b_{1,\beta},\ldots,b_{N,\beta}).$ Then

\item{{\rm(1)}}~for $m\in\{1,2,\ldots\}$ and $z\in\mathbb{D}^{n},$
\be\label{eq-10}\left\|\sum_{|\beta|=m}a_{\beta}z^{\beta}+\sum_{|\beta|=m}a_{\beta}z^{\beta}\right\|\leq\frac{4}{\pi};\ee

\item{{\rm(2)}}~ for $z\in\mathbb{D}^{n},$ \be\label{eq-11} \|f(0)\|^{2}+
\sum_{|\beta|=1}^{\infty}\left(\|a_{\beta}\|^{2}+\|b_{\beta}\|^{2}\right)\leq1.
\ee

\end{lem}

\bpf We first prove (\ref{eq-10}). Since
$$\sum_{|\beta|=m}a_{\beta}z^{\beta}=\frac{1}{2\pi}\int_{0}^{2\pi}f(e^{i\theta}z)e^{-im\theta}\,d\theta$$
and
$$\sum_{|\beta|=m}\overline{b}_{\beta}\overline{z}^{\beta}=\frac{1}{2\pi}\int_{0}^{2\pi}f(e^{i\theta}z)e^{im\theta}\,d\theta,$$
by Lemma \ref{lem-1}, we see that

\begin{eqnarray*}
\left\|\sum_{|\beta|=m}a_{\beta}z^{\beta}+\sum_{|\beta|=m}a_{\beta}z^{\beta}\right\|&=&\left\|\frac{1}{2\pi}
\int_{0}^{2\pi}f(e^{i\theta}z)(e^{-im\theta}+e^{im\theta})\,d\theta\right\|\\
&\leq&\frac{1}{2\pi}
\int_{0}^{2\pi}\left\|f(e^{i\theta}z)\right\|\left|e^{-im\theta}+e^{im\theta}\right|\,d\theta\\
&\leq&\frac{1}{\pi}\int_{0}^{2\pi}|\cos n\theta|\,d\theta\\
&=&\frac{4}{\pi}.
\end{eqnarray*}

Now we prove (\ref{eq-11}). For $\xi_{j}\in\mathbb{D},$ %and
%$\theta_{j}\in[0,2\pi]$,
let $\xi=(\xi_{1},\ldots,\xi_{n})$, where $j\in\{1,\ldots,n\}$. Then

\begin{multline*}
\frac{1}{(2\pi)^{n}}\int_{0}^{2\pi}\cdots\int_{0}^{2\pi}\left\|f(\xi_{1}e^{i\theta_{1}},\ldots,\xi_{n}e^{i\theta_{n}})\right\|^{2}d\theta_{1}\cdots
d\theta_{n}\\
 =\|f(0)\|^{2} + \sum_{|\beta|=1}^{\infty}\left(\|a_{\beta}\|^{2}+\|b_{\beta}\|^{2}\right)\big|\xi^{\beta}\big|^{2}
\leq 1.
\end{multline*} By letting $\xi\rightarrow\partial\mathbb{D}^{n}$,
we get the desired result. \epf

\subsection*{Proof of Theorem \ref{thm-2}}
 We first prove  (\ref{eq17}). For any fixed $z\in\mathbb{D}^{n}$, let
$\phi_{z}$ be a holomorphic automorphism of $\mathbb{D}^{n}$ with
$\phi_{z}(0)=z.$ For $\varsigma\in\mathbb{D}^{n},$ let
$F(\varsigma)=f(\phi_{z}(\varsigma)).$ Then
$$DF(\varsigma)=Df(\phi_{z}(\varsigma))D\phi_{z}(\varsigma)~\mbox{and}~
\overline{D}F(\varsigma)=\overline{D}f(\phi_{z}(\varsigma))\overline{D\phi_{z}(\varsigma)}.$$
Applying Lemma \ref{thm-2} (\ref{eq-10}) to $F$, we get

\be\label{eq-13}\left\|DF(0)\varsigma+\overline{D}F(0)\overline{\varsigma}\right\|=\left\|Df(z)D\phi_{z}(0)\varsigma+
\overline{D}f(z)\overline{D\phi_{z}(0)}\overline{\varsigma}\right\|\leq\frac{4}{\pi},\ee
%which implies
%$$\left|Df(z)\theta^{T}+
%\overline{D}f(z)\overline{\theta}^{T}\right|\leq\frac{4}{\pi(1-\|z\|_{\infty}^{2})},$$
where  $\varsigma$ is regarded as a column vector and

 $$D\phi_{z}(0)=\left(\begin{array}{ccccc}
\ds 1-|z_{1}|^{2}&  0  & 0 & \cdots &  0  \\[2mm]
\ds 0&   1-|z_{2}|^{2} &  0&\cdots &   0\\[2mm]
 \ds \vdots & \vdots & \vdots & \cdots & \vdots \\[2mm]
%\ds \vdots \hspace{1cm} \vdots \\[2mm]
 \ds  0& 0  &    \cdots &
1-|z_{n-1}|^{2}&   0\\[2mm]
\ds 0&    0&   \cdots &
 0&    1-|z_{n}|^{2}
\end{array}\right).
$$
By applying (\ref{eq-13}), and by  letting
$\varsigma\rightarrow\partial\mathbb{D}^{n}$, we have
$$\left\|Df(z)\theta+
\overline{D}f(z)\overline{\theta}\right\|\leq\frac{4}{\pi(1-\min_{1\leq
k\leq n}|z_{k}|^{2})}\leq\frac{4}{\pi(1-\|z\|_{\infty}^{2})},$$
where $\theta\in\mathbb{C}^{n}$ and $\|\theta\|_{\infty}=1.$

Now we  prove (\ref{eq18}). For any fixed
$z'\in\mathbb{D}^{n}\backslash\{0\}$, letting
$$F(\zeta)=\bigg\langle f\left(\frac{\zeta z'}{\|z'\|_{\infty}}\right), \frac{f(z')}{|f(z')|}\bigg\rangle~\mbox{for}~\zeta\in\mathbb{D}.$$
By using \cite[Lemma]{HZ}, we have
$$\big|F(\|z'\|_{\infty})\big|=\|f(z')\|\leq\frac{4}{\pi}\arctan\|z'\|_{\infty}.$$
The proof of the theorem is complete. \hfill $\Box$
%Next we come to  prove Theorem \ref{thm-2} (2). For any fixed
%$z\in\mathbb{B}^{n}$, let $\varphi_{z}$ be a holomorphic
%automorphism of $\mathbb{B}^{n}$ with $\varphi_{z}(0)=z.$ For
%$\zeta\in\mathbb{B}^{n},$ let $Q(\zeta)=f(\varphi_{z}(\zeta)).$ Then
%$$DQ(\zeta)=Df(\varphi_{z}(\zeta))D\varphi_{z}(\zeta)~\mbox{and}~
%\overline{D}Q(\zeta)=\overline{D}f(\varphi_{z}(\zeta))\overline{D\varphi_{z}(\zeta)}.$$
%Applying Lemma \ref{thm-2} (\ref{eq-10}) to $F$, we get

%\be\label{eq-14}
%\left|DQ(0)\zeta+\overline{D}Q(0)\overline{\zeta}\right|=\left|Df(z)D\varphi_{z}(0)\zeta+
%\overline{D}f(z)\overline{D\varphi_{z}(0)\zeta}\right|\leq\frac{4}{\pi}.\ee
%For $\zeta\in\mathbb{B}^{n}\backslash\{0\},$ by (\ref{eq-14}), we
%have \be\label{eq-15}
%\left|Df(z)\frac{D\varphi_{z}(0)\zeta}{|D\varphi_{z}(0)\zeta|}+\overline{D}f(z)\frac{\overline{D\varphi_{z}(0)\zeta}}{|D\varphi_{z}(0)\zeta|}\right|
%\leq\frac{4}{\pi|D\varphi_{z}(0)\zeta|}\leq\frac{4}{\pi(1-|z|^{2})},
%\ee which yields that
%$$\max_{\theta\in\partial\mathbb{B}^{n}}\left|Df(z)\theta+\overline{D}f(z)\overline{\theta}\right|\leq\frac{4}{\pi(1-|z|^{2})}.$$
%\hfill $\Box$

%\begin{thm}\label{thm-3}
%Let $u=(u_{1},\cdots,u_{n})$ be a pluriharmonic mapping in
%$\mathbb{B}^{n}$, where for all $j\in\{1,\cdots,n\},$ $u_{j}$ are
%real pluriharmonic from $\mathbb{B}^{n}$ into $\mathbb{R}$.

%\end{thm}

{\bf Acknowledgements:} This research was partly  supported by
National Natural Science Foundation of China (No. 11401184 and No.
11071063), the Construct Program of the Key Discipline in Hunan
Province, and the V\"ais\"al\"a Foundation of the Finnish Academy of
Science and Letters.

\normalsize

\end{document}